\documentclass[reqno]{amsart}
\usepackage{amssymb, amsmath, amsthm, amscd, mathrsfs}

\usepackage{paralist}
\usepackage{cite}
\usepackage{bigints}
\usepackage{dsfont}
\usepackage{empheq}
\usepackage{amssymb}
\usepackage{cases}
\usepackage{enumitem}
\allowdisplaybreaks
\usepackage{bbm}
\usepackage{caption}
\usepackage{booktabs}
\usepackage{float}
\usepackage[ruled,vlined,linesnumbered,noresetcount]{algorithm2e}

\usepackage{hyperref}
\numberwithin{equation}{section}

\theoremstyle{plain}
\newtheorem{theorem}{Theorem}

\newtheorem{lemma}{Lemma}

\theoremstyle{definition}
\newtheorem{definition}[theorem]{Definition}
\newtheorem{remark}{Remark}

\newcommand{\vertiii}[1]{{\left\vert\kern-0.25ex\left\vert\kern-0.25ex\left\vert #1 
    \right\vert\kern-0.25ex\right\vert\kern-0.25ex\right\vert}}

\def \d {\mathrm{d}}

\title[Controllability and Inverse Problems for PDEs] 
      {Controllability and Inverse Problems for Hyperbolic and Dispersive Equations with Dynamic Boundary Conditions}      
 
\author{S. E. Chorfi}
\author{L. Maniar}
\author{R. Morales}

\address{S. E. Chorfi, L. Maniar, Cadi Ayyad University, UCA, Faculty of Sciences Semlalia, Laboratory of Mathematics, Modeling and Automatic Systems, B.P. 2390, Marrakesh, Morocco}
\email{s.chorfi@uca.ac.ma, maniar@uca.ma}

\address{L. Maniar, The UM6P-Vanguard Center, University Mohammed VI Polytechnic, Benguerir, Morocco}
\email{Lahcen.Maniar@um6p.ma}

\address{R. Morales, Chair of Computational Mathematics, DeustoTech, University of Deusto, Avenida de las Universidades 24,
48007 Bilbao, Basque Country, Spain.}
\email{roberto.morales@deusto.es}

\makeatletter 
\@namedef{subjclassname@2020}{%
  \textup{2020} Mathematics Subject Classification}
\makeatother

 
\begin{document}
\begin{abstract}
This review examines classical and recent results on controllability and inverse problems for hyperbolic and dispersive equations with dynamic boundary conditions. We aim to illustrate the applicability of Carleman estimates to establish exact controllability of such equations and derive Lipschitz stability estimates for inverse problems of source terms and coefficients with general dynamic boundary conditions. We highlight the challenges associated with dynamic boundary conditions compared to classical static ones. Finally, we conclude with a discussion of open problems and future research directions.
\end{abstract}
\maketitle

\tableofcontents
\section{Introduction and motivation}
Hyperbolic and dispersive equations with dynamic boundary conditions model various complex physical phenomena, such as wave propagation in macroscopic or microscopic media, where the boundary behavior changes over time. This complexity can also reflect scenarios like vibrating membranes, oscillations, or fluid flows, where the boundaries interact dynamically with the medium. On the other hand, the controllability and inverse problems of these equations with dynamic boundary conditions involve understanding how we can design controls or reconstruct unknown parameters from partial measurements.

In this paper, we examine hyperbolic and dispersive equations subject to dynamic boundary conditions, focusing on aspects such as exact controllability and inverse problems. These topics are well understood in the context of static boundary conditions like Dirichlet, Neumann, Robin, mixed, and periodic conditions, etc. We mainly focus on wave and Schr\"odinger equations as prototype models. We refer to \cite{MMS17} for the null controllability of parabolic equations with dynamic conditions. The reader can also refer to the recent review \cite{CM25} for controllability and inverse problems for parabolic equations with dynamic conditions.

\subsection{Literature on hyperbolic equations}
The exact controllability problems of wave equations with static boundary conditions are fairly well-studied, see for instance \cite{Lions1988}, \cite{Zua24} and the references therein. Namely, the controllability problem
\begin{align}
    \label{cwave}
    \begin{cases}
        y_{tt} -\Delta y=0&\text{ in }\Omega\times (0,T),\\
        y=\mathbbm{1}_{\Gamma_\star} v&\text{ in }\partial \Omega\times (0,T),\\
        y(\cdot,0)=y_0, \; y_t(\cdot,0)=y_1&\text{ in }\Omega,
    \end{cases}
\end{align}
where $T>0$ is a fixed time, $\Omega\subset \mathbb{R}^n$ is a bounded domain with boundary $\partial \Omega$, $y$ is the displacement state, $y_0, y_1$ are the initial data, $v\in L^2(\partial \Omega \times (0,T))$ is a control acting only on a portion $\Gamma_\star \subset \partial \Omega$ suitably chosen, and $\mathbbm{1}_{\Gamma_\star}$ is the indicator function of $\Gamma_\star$. The problem of exact controllability of \eqref{cwave} at time $T$ concerns whether there exists a control $v \in L^2(\partial \Omega \times (0,T))$ such that we have at time $T$:
\begin{equation}
y(\cdot, T)=0 \quad \text { and } \quad y_t(\cdot, T)=0.
\end{equation}

Similarly, one can also deal with the case of a distributed control acting on a subset $\omega \subset \Omega$ (see \cite{LLions1988}):
\begin{align}
    \label{linear:static:wave:distri}
    \begin{cases}
        y_{tt} -\Delta y=\mathbbm{1}_\omega v&\text{ in }\Omega\times (0,T),\\
        y=0&\text{ in }\partial \Omega\times (0,T),\\
        y(\cdot,0)=y_0, \; y_t(\cdot,0)=y_1&\text{ in }\Omega. 
    \end{cases}
\end{align}

\begin{figure}[!h]
    \centering
    \includegraphics[width=0.7\linewidth]{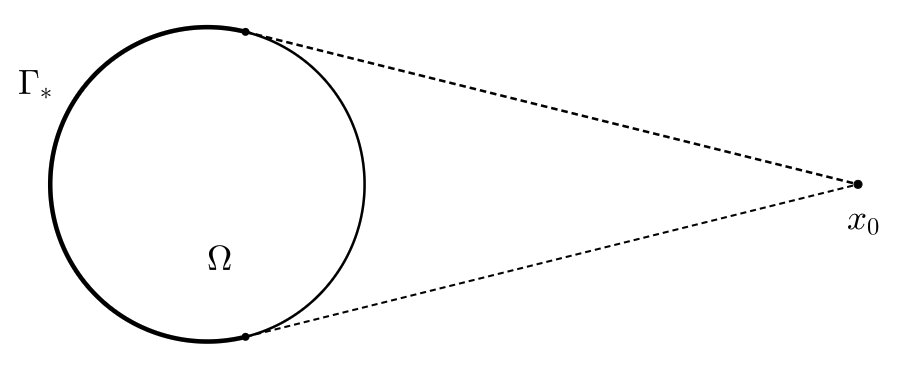}
    \caption{Domain satisfying GCC.}
    \label{fig:GCC-domain}
\end{figure}

Generally speaking, the controllability of such models depends on the so-called Geometric Control Condition (GCC in short) and a time condition; see, e.g., \cite{Bardos1992}. The GCC asserts that every ray in geometric optics encounters the boundary at a point where diffraction does not occur before time $T$. See Figure \ref{fig:GCC-domain} for a special case of GCC known as $\Gamma$-set:
\begin{equation}\label{Glions}
\Gamma_\star:=\{x\in \partial \Omega\,:\, (x-x_0)\cdot \nu(x)> 0\},
\end{equation}
where $\nu$ is the outward unit normal field to $\partial \Omega$. Intuitively, if you stand at $x_0 \notin \overline{\Omega}$ and shine a light, $\Gamma_\star$ is the shadowed part of $\partial \Omega$. If there exists $x_0 \in \mathbb{R}^n$ such that $\Gamma_\star$ is given by \eqref{Glions} and $T>2 \max\limits _{x \in \overline{\Omega}}\left|x-x_0\right|$, then it can be proved by the multiplier method that \eqref{cwave} is exactly controllable at time $T$ (see \cite{Ho86}). We refer to \cite{Os01} for an alternative approach using rotated
multipliers and providing an extension of the set $\Gamma_\star$.

Dynamical boundary conditions in the context of waves, also called nonlocally reacting kinetic conditions, originate from a kinetic energy function that includes boundary terms (potential energy may also be considered); see, e.g., \cite{H17} for more details.

The modeling and well-posedness of hyperbolic systems with dynamic boundary conditions have been the subjects of active research in the last decades. For instance, the paper \cite{Hi89} discusses evolution equations with dynamic boundary conditions and proves the existence and uniqueness of solutions using semigroup theory. The article \cite{BE04} studies the well-posedness of abstract wave equations with generalized Wentzell boundary conditions using cosine families. In \cite{FGGR05}, the authors explore the well-posedness of a non-autonomous one-dimensional wave equation with general Wentzell boundary conditions. The work \cite{GGG03} examines oscillatory boundary conditions for acoustic wave equations. Notably, the seminal paper \cite{Go06} provides a derivation and physical interpretation of general dynamic boundary conditions. The authors of \cite{GL14} investigate analyticity and Gevrey class regularity for a strongly damped wave equation with hyperbolic dynamic boundary conditions. In \cite{Gu20}, the study focuses on the well-posedness of general hyperbolic mixed problems with dynamic and Wentzell boundary conditions. The author of \cite{Mu11} investigates damped wave equations with dynamic boundary conditions and some qualitative properties of their solutions. The paper \cite{Vi17} addresses some nonlinear damped wave equations with hyperbolic dynamical boundary conditions, by studying local and global well-posedness.

Controllability and inverse problems for wave-like systems with static boundary conditions have been studied extensively; see, for instance, \cite{Bardos1992, FI96, Im02, IY01, IY012}, to cite a few. Dynamical boundary conditions have been recently explored in several papers. We summarize a few, but the list is not exhaustive. The paper \cite{AL03} discusses the exact controllability of structural acoustic interactions using microlocal estimates. In \cite{BDEM22}, the authors present a unified strategy for the observability of waves in an annulus with various boundary conditions (including dynamic conditions) based on resolvent estimates. The work \cite{BT18} establishes the observability of a wave equation with Ventcel (Wentzell) conditions by proving a variant of Mehrenberger's inequality. The paper \cite{CNR23} investigates Lyapunov functions for linear damped wave equations with dynamic boundary conditions to prove exponential stability results. 

As a closely related topic, we mention \cite{HZ95}, where the authors explore the exact controllability and stabilization of a vibrating string with an interior point mass by studying observability estimates and moment problems. The paper \cite{Wang22} investigates the controllability of a hyperbolic system with dynamic boundary conditions using the moment method. The work \cite{Wang17} focuses on local exact boundary controllability for 1-D quasilinear wave equations with dynamical boundary conditions by a constructive approach. The paper \cite{Wang19} extends the exact boundary controllability to a coupled system of quasilinear wave equations. The article \cite{Ma14} examines the controllability in polygonal and polyhedral domains using the multiplier technique. Using a similar method, the note \cite{LX16} discusses the initial boundary value problem for wave equations with dynamic boundary conditions. In \cite{VFP20}, the authors study the control of a wave equation with a dynamic boundary condition through a Neumann control. The research \cite{ZG19} is devoted to the uniform stabilization of semilinear wave equations with localized internal damping. The exponential stabilization is proved by constructing suitable Lyapunov functionals.

The authors of \cite{GT17} proved Carleman inequalities for wave equations with oscillatory boundary conditions, but the proof relies on some inaccurate calculations. More recently, the article \cite{CGMZ24} corrects these calculations and proves a refined Carleman estimate by modifying the weight function. The authors then address the Lipschitz stability for an inverse source problem following the seminal method of \cite{BK81} (see also the book \cite{BY17}). The authors also prove a new controllability result using a Dirichlet control. In addition, an inverse source problem for a one-dimensional wave equation with dynamic boundary conditions was numerically investigated in \cite{CGMZ23} using optimization techniques. For numerical analysis of wave equations with dynamic boundary conditions, we refer to the recent studies \cite{H17, HK21}.

\subsection{Literature on dispersive equations}
There is a vast literature concerning the exact controllability of the Schr\"odinger-type equations with static boundary conditions. To fix some ideas, let us consider the following (controlled) problem in a bounded domain $\Omega\subset \mathbb{R}^n$ with a boundary $\partial\Omega$:
\begin{align}
    \label{linear:static:schrodinger:bdy}
    \begin{cases}
        iy_t -\Delta y=0&\text{ in }\Omega\times (0,T),\\
        y=\mathbbm{1}_{\Gamma_\star} v&\text{ in }\partial \Omega\times (0,T),\\
        y(\cdot,0)=y_0&\text{ in }\Omega. 
    \end{cases}
\end{align}
Here, $i$ is the imaginary unit, $y:\Omega \times (0,T)\rightarrow \mathbb{C}$ is the state describing the amplitude of a quantum particle, $y_0$ is the initial condition, and the control $v\in L^2(\partial \Omega\times (0,T))$ is active only on $\Gamma_\star \subset \partial \Omega$. We can also consider the case of a distributed control acting on a subset $\omega \subset \Omega$:
\begin{align}
    \label{linear:static:schrodinger:distri}
    \begin{cases}
        iy_t -\Delta y=\mathbbm{1}_\omega v&\text{ in }\Omega\times (0,T),\\
        y=0&\text{ in }\partial \Omega\times (0,T),\\
        y(\cdot,0)=y_0&\text{ in }\Omega. 
    \end{cases}
\end{align}
We first mention the paper \cite{Leb92} where the author proved that the GCC for the exact controllability of the wave equation implies the exact controllability for the Schr\"odinger equation in any time $T>0$. In \cite{Ma90, Ma94}, the exact controllability was proved using multiplier techniques. Consequently, these results have been used in \cite{Las1992} and \cite{MaZ94} to study the problem of boundary feedback stabilization. 

It is well-known that, when GCC is not satisfied at any time $T$, one can still achieve controllability results for the Schr\"odinger equation. This can be obtained using auxiliary controllability results for the plate equation with hinged boundary conditions, see for instance \cite{Ja1988} and \cite{Bu1993}. Roughly speaking, the strategy is to use the fact that the plate operator can be written as a composition of two conjugate Schr\"odinger operators and transfer the exact controllability results from one system to another as in \cite{Leb92}. Furthermore, the paper \cite{Phung2001} estimates the cost of approximate controllability for the Schr\"odinger equation when GCC is not satisfied, and provides how the control size depends on the control time.

The exact controllability of the systems \eqref{linear:static:schrodinger:bdy} and \eqref{linear:static:schrodinger:distri} is equivalent to suitable observability inequalities for the following adjoint system (see, e.g., \cite{TW09}):
\begin{align}
    \label{adjoint:schrodinger:linear:Dirichlet}
    \begin{cases}
        i z_t-\Delta z=0&\text{ in }\Omega\times (0,T),\\
        z=0&\text{ on }\partial \Omega\times (0,T),\\
        z(\cdot,T)=z_0&\text{ in }\Omega.
    \end{cases}
\end{align}
More precisely, the exact controllability of \eqref{linear:static:schrodinger:bdy} is equivalent to the existence of a constant $C_{\text{obs}}>0$ such that the following inequality holds
\begin{align}
    \|z_0\|_{L^2(\Omega)}^2 \leq C_{\text{obs}}\int_0^T\int_{\Gamma_\star} \left|\partial_\nu z\right|^2\,\d x\,\d t \quad \forall\, z_0\in L^2(\Omega),
\end{align}
where $z$ is the associated solution of \eqref{adjoint:schrodinger:linear:Dirichlet}. In the same manner, regarding the controllability of \eqref{linear:static:schrodinger:distri}, we shall prove that there is $C_{\text{obs}}>0$ such that
\begin{align}
    \|z_0\|_{L^2(\Omega)}^2 \leq C_{\text{obs}}\int_0^T\int_{\omega} |z|^2\,\d x\,\d t \quad \forall\, z_0\in L^2(\Omega).
\end{align}

In \cite{Baudouin02}, the authors used Carleman inequalities (in the context of inverse problems) for a Schr\"odinger equation with real potentials. After that, in \cite{Las2004p1} and \cite{Las2004p2}, some Carleman inequalities and observability estimates in $H^1(\Omega)$ and $L^2(\Omega)$, respectively, are obtained for non-conservative Schr\"odinger systems with mixed boundary conditions. We also mention \cite{Liu11} and \cite{HKS24} where some inverse problems for the magnetic Schr\"odinger equation are considered. 

\begin{figure}[!h]
    \centering
    \includegraphics[width=0.6\linewidth]{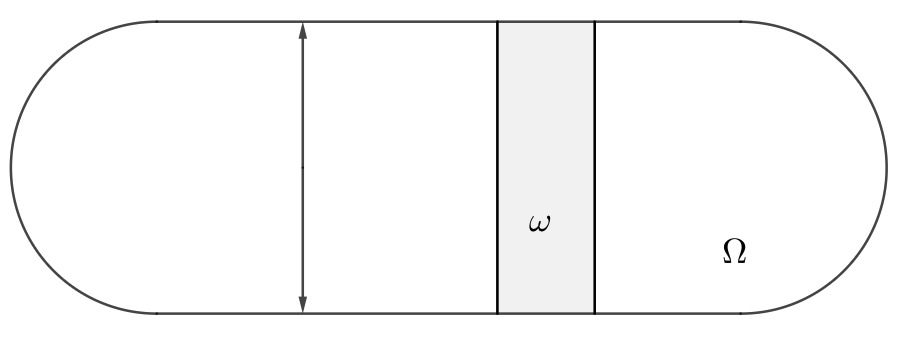}
    \caption{Domain with an interior observation region which does not satisfy GCC.}
    \label{fig:stadium:interior}
\end{figure}

\begin{figure}[!h]
    \centering
    \includegraphics[width=0.6\linewidth]{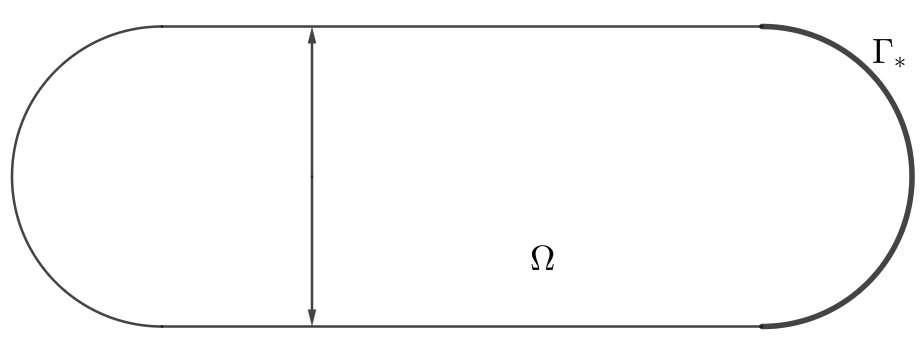}
    \caption{Domain with a boundary observation region which does not satisfy GCC.}
    \label{fig:stadium:boundary}
\end{figure}

But, as we mentioned before, the exact controllability of the Schr\"odinger equation can be obtained in domains in which the GCC condition is not satisfied. See the example of a `stadium' in Figures \ref{fig:stadium:interior} and \ref{fig:stadium:boundary} where the observation region is taken in the interior and on a part of the boundary, respectively. In both cases, a trapped ray (represented by an arrow) avoids the observation region.

In \cite{MOR08}, in the context of inverse problems, the authors obtained Carleman inequalities for the Schr\"odinger equation with Dirichlet boundary conditions, by using suitable degenerate weights. Then, with these estimates at hand, the authors solved the problem of retrieving a time-independent potential for the Schr\"odinger equation from a single boundary or internal measurement. In a similar way, in \cite{LMT16}, the exact controllability of a system of Schr\"odinger equations in cascade form is studied. Using the degenerate weights proposed in \cite{MOR08}, they achieve a Carleman estimate for the associated adjoint system. Moreover, the recent paper \cite{Ch25} establishes uniqueness results for a singular Schr\"odinger equation with an inverse square potential to study approximate controllability and inverse source problems.

Regarding the case of the Schr\"odinger equation with dynamic boundary conditions, we mention \cite{CCLL16}, where the authors studied the well-posedness and stability issues for some nonlinear Schr\"odinger equations with Wentzell boundary conditions. A complex Ginzburg-Landau equation with dynamic boundary conditions has been studied in \cite{CO18}. Moreover, in \cite{MM23}, the exact controllability of the linearized case has been analyzed under suitable geometric assumptions. In \cite{CHM24}, an inverse problem of retrieving the source terms of a Schr\"odinger equation with dynamic boundary conditions from final time measurements is considered. Recently, in \cite{CMM25}, the authors studied an inverse coefficient problem for such Schr\"odinger-type models. Then, using a suitable Carleman estimate and the Bukhgeim-Klibanov method, a Lipschitz stability estimate is obtained from observations on a part of the boundary.

\subsection{Outline of the paper}

This paper is organized as follows. In Section \ref{section:preliminaries}, we set up the main assumptions of the article and some needed function spaces. In Section \ref{sec2}, some controllability and stability results for hyperbolic equations with dynamic boundary conditions are presented. In Section \ref{sec3}, controllability results and inverse problems for the Schr\"odinger equations with dynamic boundary conditions are considered. Finally, in Section \ref{sec4}, some open problems are proposed. 

\section{Preliminaries}
\label{section:preliminaries}
\subsection{Geometric setting}\label{gsec}
Let $T>0$ be a fixed time. Throughout the paper, we assume $\Omega\subset \mathbb{R}^n$ is a bounded domain ($n$ is an integer) with a boundary $\Gamma$, $\Omega = \Omega_0 \setminus \overline \Omega_1$, where $\Omega_0$ is an open bounded set with $C^2$-boundary and $\Omega_1$ is an open strongly convex set such that $\overline{\Omega_1}\subset \Omega_0$. In particular, $\Gamma=\Gamma^0\cup \Gamma^1$, with $\Gamma^{k} = \partial \Omega_k$ for $k=0,1$ and $\Gamma^0\cap \Gamma^1=\varnothing$. By translation, we may assume $0 \in \Omega_1$ so that $0\notin \overline{\Omega}$. See Figure \ref{fig:domain}:
\begin{figure}[!h]
    \centering
    \includegraphics[width=0.6\linewidth]{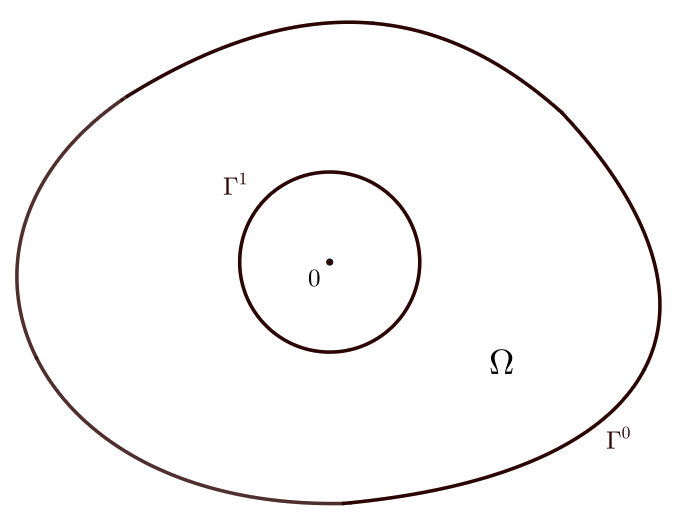}
    \caption{Domain satisfying the geometric setting.}
    \label{fig:domain}
\end{figure}

As is \cite{MM23}, we consider the Minkowski gauge functional of the set $\Omega_1$ given by
\begin{equation} \label{Mink}
\mu(x) = \inf \{\lambda \, \colon \lambda > 0 \, \text{ and } \, x \in \lambda \Omega_1\}, \qquad x \in \mathbb{R}^n.
\end{equation}
To guarantee sufficient regularity of the function $\mu$, we assume that $\Gamma^1$ belongs to the class $C^4$. Additionally, we consider
\begin{equation} \label{psi0}
\psi_{0}(x)=\mu^2(x), \qquad x \in \mathbb{R}^n.
\end{equation}
\begin{lemma}\label{lmw}
The weight function $\psi_{0}$ satisfies the following properties:
\begin{itemize}
	\item[$\mathrm{(i)}$] $\psi_{0}\in C^4(\overline{\Omega})$,
	\item[$\mathrm{(ii)}$] $\psi_{0} = 1$ on $\Gamma^1$,
	\item[$\mathrm{(iii)}$] $\nabla \psi_{0}(x) \neq 0$ for all $x\in \overline{\Omega}$,
	\item[$\mathrm{(iv)}$] There exists $\rho > 0$ such that $\nabla^2 \psi_{0} (\xi, \xi) \geq 2\rho |\xi|^2$  in $\overline{\Omega}$ for all $\xi \in \mathbb{R}^n$,
    \item[$\mathrm{(v)}$] $\partial_\nu \psi_0 < 0$ on $\Gamma^1$.
\end{itemize}
\end{lemma}
We also define the subset
\begin{align}
    \label{def:Gamma:star}
    \Gamma_\star:=\{x\in \Gamma\,:\, \partial_\nu \psi_0(x)\geq 0\}\subseteq \Gamma^0, 
\end{align}
which will be the control/measurement region throughout the paper.

\subsection{Functional setting}
We denote the Lebesgue measure on $\Omega$ by $\d x$ and the surface measure on $\Gamma$ by $\d S$. Additionally, we consider the functional space: $\mathcal{H} := L^2(\Omega,\d x) \times L^2\left(\Gamma^1, \d S\right)$. It is a real Hilbert space equipped with the inner product given by
\begin{align*}
&\langle (u,u_\Gamma),(v,v_\Gamma)\rangle_{\mathcal{H}}:=\int_\Omega u v\,\d x + \int_{\Gamma^1} u_\Gamma v_\Gamma\,\d S,
\end{align*}
which will be used for the wave equation. For the Schr\"odinger equation, the function spaces will refer to complex-valued functions, and the space $\mathcal{H}$ will be endowed with the Hermitian product 
\begin{align*}
    \left( (u,u_\Gamma),(v,v_\Gamma) \right)_{\mathcal{H}}:=\int_\Omega u\overline{v}\,\d x + \int_{\Gamma^1} u_\Gamma \overline{v_\Gamma}\,\d S.
\end{align*}
We also consider
\begin{align*}
H_{0,\Gamma^0}^k(\Omega) &:=\left\{y \in H^k(\Omega): y=0 \text { on } \Gamma^0\right\},\\
\mathcal{E}^k &:=\left\{\left(y, y_{\Gamma}\right) \in H_{ 0,\Gamma^0}^k(\Omega) \times H^k\left(\Gamma^1\right): y_{\Gamma}=y_{\mid_{\Gamma^1}}\right\}, \qquad k=1,2.
\end{align*}
We also set $\mathcal{E}^1 = \mathcal{E}$ and $\mathcal{E}^{-1}$ for its topological dual with respect to the pivot space $\mathcal{H}$.

We point out that $\mathcal{E}^1$ is a Hilbert space in $\mathbb{C}$ equipped with the inner product
\begin{align*}
    \left( (y,y_\Gamma), (z,z_\Gamma)\right)_{\mathcal{E}^1}:=\int_\Omega \nabla y\cdot \nabla \overline{z}\,\d x + \int_{\Gamma^1} \nabla_\Gamma y_\Gamma \cdot \nabla_\Gamma \overline{z_\Gamma}\,\d S
\end{align*}
for all $(y,y_\Gamma),(z,z_\Gamma)\in \mathcal{E}^1$. 
\section{Hyperbolic equations with dynamic boundary conditions}\label{sec2}
We review exact controllability and inverse problems for wave equations with dynamic boundary conditions. More details can be found in the paper \cite{CGMZ24}.

\subsection{Boundary controllability}
In this section, we provide a concise discussion on the exact boundary controllability of the following system
\begin{align}
	\label{ctrlpb}
	\begin{cases}
		\partial_{tt} y-d\Delta y + q_{\Omega}(x,t) y=0, &\text{ in } \Omega \times (0,T),\\
		\partial_{tt} y_\Gamma -\delta \Delta_\Gamma y_\Gamma + d\partial_\nu y + q_{\Gamma}(x,t) y_\Gamma=0, &\text{ on }\Gamma^1\times (0,T) ,\\
        y_\Gamma = y_{\mid_{\Gamma^1}}, &\text{ on }\Gamma^1\times (0,T),\\
		y=\mathds{1}_{\Gamma_\star} v,&\text{ on }\Gamma^0\times (0,T),\\
		(y(\cdot,0),y_\Gamma (\cdot,0))=(y_0,y_{0,\Gamma}), &\text{ in }\Omega\times \Gamma^1,\\
        (\partial_t y(\cdot,0),\partial_t y_\Gamma (\cdot,0))=(y_1,y_{1,\Gamma}), &\text{ in }\Omega\times \Gamma^1,
	\end{cases}
\end{align}
with a single control $v\in L^2(\Gamma \times (0,T))$, where $q_{\Omega} \in L^{\infty}(\Omega\times (0,T))$, $q_{\Gamma} \in L^{\infty}\left(\Gamma^1\times (0,T)\right)$, and $\Gamma_\star$ is the control region defined in \eqref{def:Gamma:star}. The function $y_{\mid_{\Gamma}}$ is the trace of $y$ on $\Gamma$, $\Delta_\Gamma$ is the Laplace-Beltrami operator, and $\partial_\nu y$ is the normal derivative with the outer unit normal vector field $\nu$.
Moreover, $\sqrt{d}$ and $\sqrt{\delta}$ represent the wave velocities in the bulk and on the boundary, respectively ($d,\delta>0$ being constant).

When the initial data are less regular, the solution of \eqref{ctrlpb} should be interpreted in the transposition sense.

Simultaneously with the system \eqref{ctrlpb}, we consider its adjoint backward system
\begin{align}
	\label{adjpb}
	\begin{cases}
		\partial_{tt} z-d\Delta z + q_{\Omega}(x,t) z=0, &\text{ in } \Omega\times (0,T),\\
		\partial_{tt} z_\Gamma -\delta \Delta_\Gamma z_\Gamma + d\partial_\nu z + q_{\Gamma}(x,t) z_\Gamma=0, &\text{ on }\Gamma^1\times (0,T) ,\\
        z_\Gamma = z_{\mid_{\Gamma^1}}, &\text{ on }\Gamma^1\times (0,T) ,\\
		z=0,&\text{ on }\Gamma^0\times (0,T),\\
		(z(\cdot,T),z_\Gamma (\cdot,T))=(z_0,z_{0,\Gamma}), &\text{ in }\Omega\times \Gamma^1,\\
        (\partial_t z(\cdot,T),\partial_t z_\Gamma (\cdot,T))=(z_1,z_{1,\Gamma}), &\text{ in }\Omega\times \Gamma^1.
	\end{cases}
\end{align}
Next, we define the minimal time
\begin{equation}\label{mt}
T_*=\dfrac{2}{\min\left(\sqrt{\rho d},\sqrt{\frac{d(\delta-d)}{8\delta}}\right)} \left(\underset{x\in\overline{\Omega}}{\max}\,\psi_0(x)-\underset{x\in\overline{\Omega}}{\min}\,\psi_0(x)\right)^{\frac{1}{2}}.
\end{equation}
Note that the time $T_*$ that will be needed for observability depends on the Minkowski gauge function $\mu$, the parameter $\rho$ in Lemma \ref{lmw}-$\mathrm{(iv)}$, the constant $d$, and the difference $\delta - d$.

By a standard duality argument, the exact controllability of system \eqref{ctrlpb} is equivalent to satisfying the observability inequality for the solution of the adjoint system \eqref{adjpb}, as stated in the following result.
\begin{theorem}\label{thmobs} 
We assume that $T>T_*$ and $\delta>d$. Then there exists a constant $C>0$ such that for all $(z_0,z_{0,\Gamma},z_1,z_{1,\Gamma})\in \mathcal{E}\times \mathcal{H}$, the solution of \eqref{adjpb} satisfies the observability inequality
\begin{align}
    \begin{split} 
    \frac{1}{2} &\int_{\Omega}\left(\left|z_1(x)\right|^2+d|\nabla z_0(x)|^2\right) \d x +\frac{1}{2} \int_{\Gamma^1}\left(\left|z_{1,\Gamma}(x)\right|^2+\delta\left|\nabla_{\Gamma} z_{0,\Gamma}(x)\right|^2\right) \d S\\ 
    \leq &C \int_{\Gamma_\star\times (0,T)} |\partial_\nu z|^2 \d S \d t.
    \end{split} 
\end{align}
\end{theorem}
Next, we set the following notations:
\begin{align}\label{Inclu}
\Omega_T = \Omega \times (-T,T) \quad &\text{ and } \quad \Gamma^k_{T} = \Gamma^k \times (-T,T) \quad \text{ for } k=0,1, \notag\\
&\Gamma_\star :=\left\{ x\in\Gamma \colon \;\partial_\nu\psi_0(x)\geq 0 \right\} \subseteq \Gamma^{0}.
\end{align}
Now, we introduce the Carleman weights for the wave system \eqref{ctrlpb}:
\begin{align}
\label{def:weight:functions}
\psi(x,t)=\psi_{0}(x)-\beta t^2+C_{1}\quad \text{and} \quad\varphi(x,t)=e^{\lambda \psi(x,t)}, \quad  (x,t)\in \overline{\Omega}\times [-T,T],
\end{align}
where $0<\beta<\rho d $, $C_{1}$ is chosen so that $\psi>1$, and $\lambda$ is a parameter (sufficiently large).

The proof of Theorem \ref{thmobs} is based on the following Carleman estimate:
\begin{theorem}\label{Thm:Carleman}
There exist positive constants $C$, $s_1$ and $\lambda_1$ such that for all $ \lambda\geq \lambda_1$ and $ s\geq s_1$, the following Carleman estimate holds
\begin{align}\label{carleman}
	&\int_{\Omega_T} {e^{2s\varphi} \left(s^3\lambda^3\varphi^3|y|^2 +s\lambda \varphi |\nabla y|^2 +s\lambda \varphi |\partial_t y|^2 \right) }\mathrm{d}x \mathrm{d}t \notag\\
	&+\int_{\Gamma_{T}^{1}} e^{2s\varphi} \left(s^3\lambda^3\varphi^3 |y_\Gamma|^2 + s\lambda \varphi|\partial_\nu y|^2 + s\lambda \varphi|\partial_t y_\Gamma|^2 \right. \notag\\
    & \hspace{2cm} + \left. [d(\delta-d)-8\beta\delta]s\lambda \varphi |\nabla_\Gamma y_\Gamma|^2 \right)\mathrm{d}S\mathrm{d}t \notag\\
	\leq & C\int_{\Omega_{T}}{e^{2s\varphi}|f|^2}\mathrm{d}x\mathrm{d}t
	+ C \int_{\Gamma_{T}^{1}}{e^{2s\varphi} |g|^2}\mathrm{d}S\mathrm{d}t+C s\lambda \int_{-T}^T\int_{\Gamma_\star}{e^{2s\varphi} \varphi|\partial_\nu y|^2}\mathrm{d}S\mathrm{d}t \\
	& +  \mathcal{I}_{\pm T} \notag
\end{align}   
where
$(y,y_\Gamma)\in H^2\left(-T,T; \mathcal{H}\right)\cap L^2(-T,T;\mathcal{E})$ satisfies
$$
f :=\partial_{tt} y-d\Delta y  +q_{\Omega}y\in L^2(\Omega_T), \quad\quad  g :=\partial_{tt} y_{\Gamma} - \delta \Delta_\Gamma y_\Gamma +d\partial_\nu y + q_{\Gamma} y_\Gamma\in L^2\left(\Gamma_T^1\right),
$$
and
\begin{align}
\mathcal{I}_{\pm T} &= Cs\lambda\int_{\Omega}\mathrm{e}^{2s\varphi(x,T)}\varphi(x,T)(|\partial_t y(x,T)|^2+|\nabla y(x,T)|^2)\mathrm{d}x\notag\\
&\hspace{-0.7cm} +Cs^3\lambda^3\int_{\Omega}e^{2s\varphi(x,T)}\varphi^3(x,T)|y(x,T)|^2\mathrm{d}x \notag\\
    & \hspace{-0.7cm} +Cs\lambda\int_{\Omega}e^{2s\varphi(x,-T)}\varphi(x,-T) \left(|\partial_t y(x,-T)|^2+|\nabla y(x,-T)|^2\right)\mathrm{d}x \notag\\
    & \hspace{-0.7cm} +C s^3\lambda^3\int_{\Omega}e^{2s\varphi(x,-T)}\varphi^3(x,-T)|y(x,-T)|^2\mathrm{d}x+Cs\lambda \int_{\Gamma^1}e^{2s\varphi(x,T)}\varphi(x,T)|\partial_t y_\Gamma(x,T)|^2\mathrm{d}S \notag\\
    & \hspace{-0.7cm} +Cs\lambda \int_{\Gamma^1}e^{2s\varphi(x,-T)}\varphi(x,-T)|\partial_t y_\Gamma(x,-T)|^2\mathrm{d}S+Cs^3\lambda^3\int_{\Gamma^1}e^{2s\varphi(x,T)}\varphi^3(x,T)|y_\Gamma(x,T)|^2\mathrm{d}S \notag\\
    & \hspace{-0.7cm} +Cs^3\lambda^3\int_{\Gamma^1}e^{2s\varphi(x,-T)}\varphi^3(x,-T)|y_\Gamma(x,-T)|^2\mathrm{d}S.\notag
\end{align}
\end{theorem}

\begin{remark}
Theorem \ref{thmobs} extends and improves \cite[Theorem 2.2]{BDEM22} by giving a sharp lower bound for the time required for observability (see \cite[Section 5.2]{BDEM22}). For instance, if $\Omega_1=B_{1}$ (the unit open ball), then we can choose $\rho=1$. Hence, we obtain the observability with an explicit minimal time
$$T_*=\frac{2}{\min\left(\sqrt{d},\sqrt{\frac{d(\delta-d)}{8\delta}}\right)}\left(\underset{x\in\overline{\Omega}}{\max}\,|x|^2-\underset{x\in\overline{\Omega}}{\min}\,|x|^2\right)^{\frac{1}{2}}.$$
\end{remark}

\begin{remark}\label{rmkobs}
In the case when $\delta<d$, it has been shown in \cite[Theorem 2.4]{BDEM22} that the corresponding observability inequality fails at any time $T > 0$. The assumption $\delta>d$ has also been considered in \cite{MM23} for the controllability of a Schr\"odinger equation with a dynamic boundary condition. See Section \ref{subsch} and open problem 6.1 for further discussion.
\end{remark}

The main result of this section reads as follows.
\begin{theorem}
We assume that $T>T_*$ and $\delta>d$. Then the system \eqref{ctrlpb} is exactly controllable in time $T$ with a control acting on $\Gamma_\star$, i.e., for all initial data $(y_0,y_{0,\Gamma},y_1,y_{1,\Gamma})\in \mathcal{H} \times \mathcal{E}^{-1}$, there exists a control function $v\in L^2(\Gamma_\star\times (0,T))$ such that the solution of system \eqref{ctrlpb} satisfies
$$(y(T), y_\Gamma(T))=(0,0) \text{ and } \left(\partial_t y(T), \partial_t y_{\Gamma}(T)\right)=(0,0) \text{ in } \Omega \times \Gamma^1.$$
\end{theorem}
The controllability result presented above is optimal as it relies on a single control force applied to the Dirichlet boundary. This contrasts with \cite[Section 4]{GT17}, where two controls are employed. A similar observation has been made in \cite{BT18} for the case of a multi-dimensional rectangle. Additionally, we refer to \cite{Ma14} for results concerning domains where $\Omega$ is a polygon or polyhedron.
\subsection{Inverse problems}
In this section, we review recent results on inverse source problems of wave systems with dynamic boundary conditions. Consider the following wave system
\begin{align}
	\label{intro:problem:01}
	\begin{cases}
		\partial_{tt} y-d\Delta y + q_{\Omega}(x) y=f(x,t), &\text{ in } \Omega\times(0,T),\\
		\partial_{tt} y_\Gamma -\delta \Delta_\Gamma y_\Gamma + d\partial_\nu y + q_{\Gamma}(x) y_\Gamma=g(x,t), &\text{ on }\Gamma^1\times(0,T) ,\\
        y_\Gamma = y_{\mid_{\Gamma^1}}, &\text{ on }\Gamma^1\times(0,T) ,\\
		y=0,&\text{ on }\Gamma^0\times(0,T),\\
		(y(\cdot,0),y_\Gamma (\cdot,0))=(0,0), &\text{ in }\Omega\times \Gamma^1,\\
        (\partial_t y(\cdot,0),\partial_t y_\Gamma (\cdot,0))=(0,0), &\text{ in }\Omega\times \Gamma^1.
	\end{cases}
\end{align}
\noindent{\bf Inverse Source Problem (ISP):} Can one determine, in a stable manner, the unknown source terms $(f,g)$ from the knowledge of a single measurement $\partial_\nu y_{|\Gamma_\star \times (0,T)}$, where $(y,y_\Gamma)$ is the corresponding solution of \eqref{intro:problem:01}?
Indeed, uniqueness fails in general, and one needs to restrict the set of unknown source terms. Thus, we introduce the set of admissible source terms
\begin{align}\label{eqas}
\mathcal{S}(C_0) :=\left\{(f,g)\in H^1\left(0, T ; \mathcal{H}\right): \begin{array}{ll}
\hspace{-0.2cm}|\partial_t f(x,t)| \leq C_0 |f(x,0)|, &\hspace{-0.4cm}\text{ a.e. }(x,t)\in \Omega\times (0,T) \hspace{-0.2cm}\\
\hspace{-0.2cm}|\partial_t g(x,t)| \leq C_0 |g(x,0)|, &\hspace{-0.4cm}\text{ a.e. }(x,t)\in \Gamma^1\times (0,T) \hspace{-0.2cm}
\end{array}\right\},
\end{align}
where $C_0 > 0$ is a fixed constant. The main result is the global Lipschitz stability for the inverse source problem.
\begin{theorem}\label{thm:stab}
Let $T>\frac{T_*}{2}$ and $C_0>0$. We assume that
\begin{align}
        \label{condition:gamma:d}
        \delta>d.
    \end{align}
Then there exists a positive constant $C=C(\Omega, T, C_0,\| q_{\Omega} \|_{\infty},\|q_{\Gamma}\|_{\infty})$ such that for any source term $(f,g)\in \mathcal{S}(C_0)$, we have
    \begin{equation}
        \| f\|_{L^2(\Omega\times (0,T))}+\| g\|_{L^2(\Gamma^1\times (0,T))}\leq C\|\partial_t \partial_\nu y\|_{L^2(\Gamma_\star\times (0,T))}
    \end{equation}
for any sufficiently smooth solution $(y,y_\Gamma)$ of the system \eqref{intro:problem:01}.
\end{theorem}
The proof of Theorem \ref{thm:stab} relies on the so-called Bukhgeim-Klibanov method \cite{BK81}. This is done thanks to the Carleman estimate (Theorem \ref{Thm:Carleman}) and an extension of a recent argument avoiding cut-off functions \cite{HIY20}.

One can also consider\\
\noindent{\bf Coefficient Inverse Problem (CIP):} Determine $(q,q_\Gamma)\in L^\infty(\Omega)\times L^\infty(\Gamma^1)$ from the knowledge of $\partial_\nu y$ on $\Gamma_\star \times (0,T)$, $(\Gamma_\star \subseteq \Gamma^0)$, where $(y,y_\Gamma)$ is the solution of \eqref{intro:problem:01} associated to $(q,q_\Gamma)$.

We can prove Lipschitz stability for the CIP, but we omit the proof as it is similar to the arguments of Section \ref{Sivpb}.

\section{Dispersive equations with dynamic boundary conditions} \label{sec3}
In this section, we review the main results concerning the controllability of the Schr\"odinger equation with dynamic boundary conditions. The details can be found in \cite{MM23}.

\subsection{Boundary controllability}\label{subsch}
Consider the following non-conservative Schr\"odinger equation with dynamic boundary conditions
\begin{align}
    \label{eq:schrodinger:dbc}
    \begin{cases}
        i\partial_t u +d\Delta u-\vec{p}_1\cdot \nabla u +p_0 u=0&\text{ in }\Omega\times (0,T),\\
        i\partial_t u_\Gamma+\delta \Delta_\Gamma u_\Gamma-d\partial_\nu u -\vec{p}_{\Gamma,1}\cdot \nabla_\Gamma u_\Gamma +p_{\Gamma,0} u_\Gamma =0&\text{ on }\Gamma^1\times (0,T),\\
        u_\Gamma=u_{\mid_{\Gamma^1}} &\text{ on }\Gamma^1\times (0,T)\\
        u =\mathbbm{1}_{\Gamma_\star} h &\text{ on }\Gamma^0\times (0,T),\\
        (u(\cdot,0),u_\Gamma(\cdot,0))=(u_0,u_{\Gamma,0})&\text{ in }\Omega\times \Gamma^1.
    \end{cases}
\end{align}

Here, the pair $(u,u_\Gamma)$ is the state, $h\in L^2(\Gamma\times (0,T))$ is the control acting on a portion of the boundary $\Gamma_\star\subseteq \Gamma^0$. Moreover, $d>0$ and $\delta>0$ represent the quantum diffusion in the bulk and on the boundary, respectively.

We mention that the boundary control $h$ acts only on a subset of $\Gamma^0\subset \partial \Omega$. This implies that the first equation satisfied in the bulk is controlled by $h$, while the side condition $u_\Gamma=u_{\mid_{\Gamma^1}}$ is satisfied on $\Gamma^1\times (0,T)$.  

The exact controllability of the system  \eqref{eq:schrodinger:dbc} is understood in the following way:
\begin{definition}
    We say that the system \eqref{eq:schrodinger:dbc} is exactly controllable in a functional space $\mathcal{X}$ if for every time $T>0$ and all states $(u_0,u_{\Gamma,0}), (u_1,u_{\Gamma,1})\in \mathcal{X}$, there exists a boundary control $h\in L^2(\Gamma \times (0,T))$ such that the associated state $(u,u_\Gamma)$ of \eqref{eq:schrodinger:dbc} fulfills
    \begin{align*}
        (u(\cdot,T),u_{\Gamma}(\cdot,T))=(u_1,u_{\Gamma,1})\qquad \text{ in }\Omega\times \Gamma^1. 
    \end{align*}
\end{definition}

Now, define the function $\psi_0$ as in \eqref{psi0}. For $\lambda>0$ we consider the weight functions
\begin{align}
    \label{def:weights:control:schrodinger}
    \theta(x,t):=\dfrac{e^{\lambda \psi_0(x)}}{t(T-t)},\quad \alpha(x,t):= \dfrac{K - e^{\lambda \psi_0(x)}}{t(T-t)},\qquad (x,t)\in \overline{\Omega}\times (0,T),
\end{align}
where $K> \|e^{\lambda \psi_0}\|_{L^\infty(\Omega)}$ is a constant.

Let us define the following differential operators:
\begin{align*}
    L(w)&:=i\partial_t w+d\Delta w +\vec{q}_1 \cdot \nabla w +q_0 w,\\
    L_\Gamma(w,w_\Gamma)&:=i\partial_t w_\Gamma +\delta \Delta_\Gamma w_\Gamma -d\partial_\nu w + \vec{q}_{\Gamma,1}\cdot \nabla_\Gamma w_\Gamma +q_{\Gamma,0}w_\Gamma. 
\end{align*}
The proof of exact controllability draws on the following Carleman estimate.
\begin{theorem}
Consider the drift coefficients $$(\vec{q}_1,\vec{q}_{\Gamma,1})\in [L^\infty(\Omega\times (0,T))]^n\times [L^\infty(\Gamma^1\times (0,T))]^n$$ and the potentials
    $$(q_0,q_{\Gamma,0})\in L^\infty(\Omega\times (0,T))\times L^\infty(\Gamma^1\times (0,T)).$$
    In addition, suppose that
    \begin{align}\label{delta d}
        \delta>d.
    \end{align}
    Then, there exist positive constants $C_1,\lambda_1$ and $s_1$ such that the following inequality holds
    \begin{align}
        \label{Carleman:schrodinger:dbc}
        \begin{split}
            &\int_0^T\int_\Omega e^{-2s\alpha}(s^3\lambda^4 \theta^3 |w|^2 + s\lambda \theta |\nabla w|^2 + s\lambda^2 \theta |\nabla \psi_0 \cdot \nabla w|^2)\,\d x\,\d t \\
            &+\int_0^T\int_{\Gamma^1} e^{-2s\alpha}(s^3\lambda^3 \theta^3 |w_\Gamma|^2 + s\lambda \theta(|\partial_\nu w|^2 + |\nabla_\Gamma w_\Gamma|^2) )\,\d S\,\d t\\
            \leq & C_1 \int_0^T\int_\Omega e^{-2s\alpha}|L(w)|^2\,\d x\,\d t + C_1\int_0^T\int_{\Gamma^1} e^{-2s\alpha} |L(w,w_\Gamma)|^2\,\d S\,\d t \\
            &+C_1s\lambda \int_0^T\int_{\Gamma_\star} e^{-2s\alpha} \theta|\partial_\nu w|^2\,\d S\,\d t
        \end{split}
    \end{align}
    for all $\lambda \geq\lambda_1$, $s\geq s_1$, and $(w,w_\Gamma)\in L^2(0,T; \mathcal{E}^1)$ such that
    \begin{align*}
        L(w)\in L^2(\Omega\times (0,T)),\quad L(w,w_\Gamma)\in L^2(\Gamma^1\times (0,T)),\quad \partial_\nu w\in L^2(\Gamma\times (0,T)).
    \end{align*}
\end{theorem}

Now suppose that 
\begin{align}
    \label{assumptions:p1:pg1}
    \vec{p}_1:=d\nabla \pi -i\vec{r},\qquad \vec{p}_{\Gamma,1}:=\delta \nabla_\Gamma \pi_\Gamma -i\vec{r}_{\Gamma},
\end{align}
where 
\begin{align*}
    \pi \in L^\infty(0,T;W^{3,\infty}(\Omega;\mathbb{R})) \cap W^{1,\infty}(0,T;W^{1,\infty}(\Omega;\mathbb{R})),\\
    \pi_\Gamma \in L^\infty(0,T; W^{3,\infty}(\Gamma^1;\mathbb{R}))\cap W^{1,\infty}(0,T;W^{1,\infty}(\Gamma^1;\mathbb{R})),
\end{align*}
with $\pi_{\mid_{\Gamma^1}}=\pi_\Gamma$ and 
\begin{align*}
    \vec{r}\in [L^\infty(0,T;W^{2,\infty}(\Omega;\mathbb{R}))]^n,\qquad \vec{r}_\Gamma \in [L^\infty(0,T;W^{2,\infty}(\Gamma^1;\mathbb{R}))]^n,
\end{align*}
with $\vec{r}_{\mid_{\Gamma^1}}=\vec{r}_\Gamma$ and  $\vec{r}\cdot \nu \leq 0$ on $\Gamma\times (0,T)$. We also assume that 
\begin{align}
    \label{assumptions:p0:pg0}
    (p_0,p_{\Gamma,0})\in L^\infty(0,T;W^{1,\infty}(\Omega)) \times L^\infty(0,T;W^{1,\infty}(\Gamma^1)).
\end{align}

\begin{theorem}
    \label{theorem:exact:controllability:boundary}
    Let $(\vec{p}_1,\vec{p}_{\Gamma,1})$ and $(p_0,p_{\Gamma,0})$ satisfy \eqref{assumptions:p1:pg1} and \eqref{assumptions:p0:pg0}, respectively. We assume \eqref{delta d} (i.e., $\delta > d$). Then, for all states $(u_0,u_{\Gamma,0}), (u_1,u_{\Gamma,1})\in \mathcal{E}^{-1}$, there exists a control $h\in L^2(\Gamma \times (0,T))$ such that the associated solution $(u,u_\Gamma)$ (in the sense of transposition) of \eqref{eq:schrodinger:dbc} satisfies
    \begin{align}
        (u(\cdot,T),u_\Gamma(\cdot,T))=(u_1,u_{\Gamma,1})\qquad \text{ in }\Omega\times \Gamma^1.
    \end{align}
\end{theorem}
By duality, the exact controllability of \eqref{eq:schrodinger:dbc} is equivalent to proving a suitable observability inequality for the adjoint system 
\begin{align}
    \label{adjoint:system:schr}
    \begin{cases}
        i\partial_t z +d\Delta z -\vec{p}_1\cdot \nabla z +p z=0&\text{ in }\Omega\times (0,T),\\
        i\partial_t z_\Gamma +\delta \Delta_\Gamma z_\Gamma-d\partial_\nu z_\Gamma -\vec{p}_{\Gamma,1}\cdot \nabla_\Gamma z_\Gamma +p_\Gamma z_\Gamma=0&\text{ on }\Gamma^1\times (0,T),\\
        z_\Gamma=z_{\mid_{\Gamma^1}} &\text{ on }\Gamma^1\times (0,T),\\
        z =0&\text{ on }\Gamma^0\times (0,T),\\
        (z(\cdot,T),z_\Gamma(\cdot,T))=(z_T,z_{\Gamma,T})&\text{ in }\Omega\times \Gamma^1,
    \end{cases}
\end{align}
where $(p,p_\Gamma)$ are given by 
\begin{align}
    p:=-d\Delta \pi +i\text{div}(\vec{r})+\overline{p}_0,\quad p_\Gamma:=-\delta \Delta_\Gamma \pi_\Gamma + i\text{div}_\Gamma (\vec{r}_\Gamma) -i\vec{r}\cdot \nu +\overline{p}_{\Gamma,0}.
\end{align}

More precisely, the exact controllability of \eqref{eq:schrodinger:dbc} is equivalent to proving the existence of a constant $C_{\text{obs}}>0$ such that for all $(z_{T},z_{\Gamma,T})\in \mathcal{E}^1$, the following observability inequality holds
\begin{align}\label{sobs}
    \|(z_T,z_{\Gamma,T})\|_{\mathcal{E}^1}^2 \leq C_{\text{obs}} \int_0^T\int_{\Gamma_\star} |\partial_\nu z|^2\,\d S\,\d t, 
\end{align}
where $(z,z_\Gamma)$ is the associated solution of \eqref{adjoint:system:schr}. Note that the observability inequality \eqref{sobs} is a direct consequence of the Carleman estimate \eqref{Carleman:schrodinger:dbc}.



\subsection{Inverse problems} \label{Sivpb}
First, we set the differential operators
\begin{align*}
    L(y)&:=i\partial_t y+d\Delta y +\vec{q}_1(x) \cdot \nabla y,\\
    L_\Gamma(y,y_\Gamma)&:=i\partial_t y_\Gamma +\delta \Delta_\Gamma y_\Gamma -d\partial_\nu y + \vec{q}_{\Gamma,1}(x)\cdot \nabla_\Gamma y_\Gamma. 
\end{align*}

We consider the following Schrödinger equation with dynamic boundary conditions:
\begin{align}
    \label{eq:ip:schrodinger}
    \begin{cases}
        L(y)+p(x)y=g&\text{ in }\Omega\times (0,T),\\
        L_\Gamma(y,y_\Gamma) + p_\Gamma(x)y_\Gamma =g_\Gamma &\text{ on }\Gamma^1\times (0,T),\\
        y_\Gamma=y_{\mid_{\Gamma^1}}&\text{ on }\Gamma^1\times (0,T),\\
        y =0&\text{ on }\Gamma^0\times (0,T),\\
        (y(\cdot,0),y_\Gamma(\cdot,0))=(y_0,y_{\Gamma,0})&\text{ in }\Omega\times \Gamma^1,
    \end{cases}
\end{align}
where $(g,g_\Gamma)$ is a source term, $p,p_\Gamma$ are unknown real potentials, and $(y_0,y_{\Gamma,0})$ is the pair of initial data. 

\noindent{\bf Coefficient Inverse Problem (CIP):} Can we determine the unknown potentials $$(p,p_\Gamma)\in L^\infty(\Omega;\mathbb{R})\times L^\infty(\Gamma^1;\mathbb{R})$$ from a single measurement of the flux $\partial_\nu y$ on $\Gamma_\star \times (0,T)$, $(\Gamma_\star \subseteq \Gamma^0)$, where $(y,y_\Gamma)$ is the solution of \eqref{eq:ip:schrodinger} associated with $(p,p_\Gamma)$?

To emphasize the dependence of the solution on $(p,p_\Gamma)$ in \eqref{eq:ip:schrodinger}, sometimes we shall write 
\begin{align*}
    y=y[p,p_\Gamma]\quad\text{ and }\quad y_\Gamma=y_\Gamma[p,p_\Gamma].
\end{align*}
For $m>0$ and $X\subset \mathbb{R}^n$, we introduce the set
\begin{align*}
    L_{\leq m}^\infty(X;\mathbb{R}):=\{p\in L^\infty(X;\mathbb{R})\,:\, \|p\|_{L^\infty(X)}\leq m\},
\end{align*}
and define the set of admissible potentials as
\begin{align*}
    \mathbb{L}_{\leq m}^\infty:=L_{\leq m}^\infty(\Omega;\mathbb{R})\times L_{\leq m}^\infty(\Gamma^1;\mathbb{R}).
\end{align*}
We also set $\mathbb{L}^\infty :=L^\infty(\Omega;\mathbb{C})\times L^\infty(\Gamma^1;\mathbb{C})$.

The next result establishes the Lipschitz stability for the coefficient inverse problem {\bf (CIP)}:
\begin{theorem}
    \label{thm:Lipschitz:stability:schrodinger}
    For a given $m>0$ we consider $(q,q_\Gamma)\in \mathbb{L}_{\leq m}^\infty$. Assume that $\delta$ and $d$ fulfill the condition \eqref{delta d} (i.e., $\delta >d$). In addition, let us assume that 
    \begin{align*}
        (y[q,q_\Gamma],y_\Gamma[q,q_\Gamma])\in H^1(0,T;\mathbb{L}^\infty)
    \end{align*}
    and the initial data $(y_0,y_{\Gamma,0})$ are real-valued functions satisfying 
    \begin{align*}
        |y_0|\geq r_0\text{ a.e. in }\Omega\quad \text{ and }\quad |y_{\Gamma,0}|\geq r_0\text{ a.e. on }\Gamma^1 
    \end{align*}
    for some constant $r_0>0$. Then, there exists a constant $C=C(\Omega,T,y_0,y_{\Gamma,0})>0$ such that if 
    \begin{align*}
        (\partial_\nu y[q,q_\Gamma]-\partial_\nu y[p,p_\Gamma])\in H^1(0,T;L^2(\Gamma_\star)),
    \end{align*}
    the following double inequality holds
    \begin{align}
        \dfrac{1}{C} \|(q,q_\Gamma)-(p,p_\Gamma)\|_{\mathcal{H}}\leq \|\partial_\nu y[q,q_\Gamma] - \partial_\nu y[p,p_\Gamma]\|_{H^1(0,T;L^2(\Gamma_\star))}\leq C\|(q,q_\Gamma)- (p,p_\Gamma)\|_{\mathcal{H}}
    \end{align}
    for all $(p,p_\Gamma)\in \mathbb{L}_{\leq m}^\infty$. 
\end{theorem}
Theorem \ref{thm:Lipschitz:stability:schrodinger} is also proved by adapting the Bukhgeim-Klibanov method \cite{BK81}. One of the main ingredients of this strategy is a suitable Carleman estimate for the Schr\"odinger operator with dynamic boundary conditions. 

Consider the function $\psi_0$ defined in \eqref{psi0}. We keep the same notation for the following modified weight functions 
\begin{align*}
    \theta(x,t):=\dfrac{e^{\lambda \psi_0(x)}}{(T+t)(T-t)},\quad \alpha(x,t):=\dfrac{K -e^{\lambda \psi_0(x)}}{(T+t)(T-t)}\quad \forall\, (x,t)\in \overline{\Omega}\times (-T,T),
\end{align*}
where $K > \|e^{\lambda \psi_0}\|_{L^\infty(\Omega)}$ is constant.

We introduce the differential operators
\begin{align*}
&\tilde{L}(w):=i \partial_t w+d \Delta w-\vec{\rho}_1 \cdot \nabla w+\rho_0 w \quad \text { in } \Omega_T,\\
&\tilde{L}\left(w, w_{\Gamma}\right):=i \partial_t w_{\Gamma}+\delta\Delta_{\Gamma} w_{\Gamma}-d \partial_\nu w-\vec{\rho}_{\Gamma, 1} \cdot \nabla_{\Gamma} w_{\Gamma}+\rho_{\Gamma, 0} w_{\Gamma} \quad \text { on } \Gamma^1_T.
\end{align*}

\begin{theorem}
    \label{thm:Carleman:inv:prob:Schrodinger}
    Consider $(\vec{\rho}_1,\vec{\rho}_{\Gamma,1})\in L^\infty(-T,T;[\mathbb{L}^\infty]^n)$ and $(\rho_0,\rho_{\Gamma,0})\in L^\infty(-T,T;\mathbb{L}^\infty)$. Suppose that $d,\delta>0$ satisfy the condition \eqref{delta d}. Then, there exist positive constants $C_1$, $\lambda_1$ and $s_1$ such that for all $\lambda \geq \lambda_1$ and $s\geq s_1$, we have
    \begin{align}
        \begin{split}
            &s^{3/2}\lambda^{3/2}\left(\int_\Omega e^{-2s\alpha(\cdot,0)} |w(\cdot,0)|^2\,\d x +\int_{\Gamma^1} e^{-2s\alpha(\cdot,0)} |w_\Gamma (\cdot,0)|^2\,\d S \right)\\
            \leq & C_1 \int_{\Omega_T} e^{-2s\alpha} |\tilde{L}(w)|^2\,\d x\,\d t + C_1 \int_{\Gamma^1_T} e^{-2s\alpha} |\tilde{L}(w,w_\Gamma)|^2\,\d S\,\d t \\
            &+C_1 s\lambda \int_{-T}^T\int_{\Gamma_\star} e^{-2s\alpha} \theta |\partial_\nu w|^2\,\d S\,\d t
        \end{split}
    \end{align}
    for all $(w,w_\Gamma)\in \mathcal{W}$, where
    \begin{align*}
        \mathcal{W} :=\big\{(w,w_\Gamma)\in L^2(-T,T;\mathcal{E}^1)\,:\, &(\tilde{L}(w),\tilde{L}_\Gamma (w,w_\Gamma)) \in L^2(-T,T;\mathcal{H})\\
         & \text{ and } \partial_\nu w\in L^2(-T,T;L^2(\Gamma_\star)) \big\}.
    \end{align*}
\end{theorem}

We point out that, thanks to Theorem \ref{thm:Carleman:inv:prob:Schrodinger}, one can develop the so-called Carleman-based Reconstruction algorithm (CbRec in short) to find an iterative reconstruction formula for the real potentials $(p,p_\Gamma)$ in \eqref{eq:ip:schrodinger}.

\section{Conclusion and open problems}\label{sec4}
This review summarizes the recent results concerning controllability and inverse problems for hyperbolic and dispersive PDEs with dynamic boundary conditions. The results have been obtained by proving new Carleman estimates for the corresponding operators with dynamic boundary conditions. These estimates have a particular interest since they can be used to deduce new controllability or stability results for systems of PDEs with dynamic boundary conditions.  

However, several interesting questions are still open and deserve further investigation in the future. In the following, we describe some of them.

\subsection{The case $\delta \leq d$}
As we have seen in both the wave and Schr\"odinger equations, the Carleman estimates \eqref{carleman} and \eqref{Carleman:schrodinger:dbc} have been obtained under the assumption $\delta >d$ on the wave speeds in $\Omega$ and on $\Gamma^1$. In both cases, this assumption is imposed to control the global term
\begin{align*}
    \int_{\Gamma^1_T} |\nabla_\Gamma z_\Gamma|^2\,\d S\,\d t  
\end{align*}
and put it on the left-hand side of the Carleman estimate. Thus, at this point, one might think that this assumption is technical for both equations. However, it was shown in \cite[Theorem 2.4]{BDEM22} that, in the context of the wave equation with dynamic boundary conditions posed in an annulus, the system is not controllable at any time $T > 0$ when $\delta < d$. To the best of the author's knowledge, the exact controllability of both equations remains open when $\delta = d$ (identical wave speeds in $\Omega$ and $\Gamma$). As for the Schr\"odinger equation with dynamic boundary conditions, the exact controllability in the case $\delta < d$ is still open in contrast to the wave equation. Finally, we emphasize that the stability of inverse problems for both equations also remains open when $\delta \leq d$.

\subsection{Geometric assumptions}
The results presented in Sections \ref{sec2} and \ref{sec3} are valid under the assumption that the domain $\Omega$ satisfies the geometric configuration described in Section \ref{gsec}. This assumption is a particular case of the Geometric Control Condition (GCC) and is connected to the derivation of Carleman inequalities for the wave or Schr\"odinger operators with dynamic boundary conditions. In particular, it allows us to apply the surface divergence theorem to control global terms involving tangential derivatives. Whether this geometric assumption can be relaxed while still achieving controllability or stability results remains an open question.

For the wave equation with dynamic boundary conditions, we have obtained the formula \eqref{mt} of a minimal time $T_*$ needed for exact controllability and stability. However, we do not know if this time is sharp and how it relates to the GCC time.

Another interesting question regarding the Schr\"odinger equation with dynamic boundary conditions is the identification of geometric conditions under which controllability results can be obtained in domains that do not satisfy the GCC. In particular, it would be worthwhile to establish Carleman inequalities for these models by employing suitable degenerate weights, as done in \cite{MOR08}.

\subsection{Controllability results for nonlinear models}
The results reviewed in Sections \ref{sec2} and \ref{sec3} pertain to linear models. However, as shown in the case of the heat or Ginzburg–Landau equations with dynamic boundary conditions (see, for instance, \cite{CMM23}), one can expect similar controllability results for nonlinear models, at least for semilinear wave and Schrödinger equations with dynamic boundary conditions.   
\section*{Acknowledgments}
R. Morales has received funding from the European Research Council (ERC) under the European Union’s Horizon 2030 research and innovation program (grant agreement NO: 101096251-CoDeFeL).

\end{document}